\documentclass[12pt]{article}

\usepackage{amsmath}
\usepackage{setspace}
\usepackage{array}
\usepackage{times}
\usepackage{bm}
\usepackage{graphicx}
\usepackage{latexsym}
\usepackage{color}

\usepackage{epsfig}

\usepackage[dvips, letterpaper, nohead, top=1.3in, bottom=1.3in, left=1.3in, right=1.3in]{geometry}

\newtheorem{theorem}{Theorem}
\newtheorem{lemma}{Lemma}
\newtheorem{definition}{Definition}
\newtheorem{proposition}{Proposition}
\newtheorem{corollary}{Corollary}
\newtheorem{remark}{Remark}

\newtheorem{assumption}{Assumption}
\newtheorem{procedure}{Procedure}
\def\qed{\rule{2mm}{2mm}}

\normalsize

\begin{document}

\small\normalsize
\title
{ Control of Directional Errors  in Fixed Sequence Multiple Testing}
\author
{Anjana Grandhi \\
BARDS, Merck Research Laboratories \\
Rahway, NJ 07065  \\
\and
 Wenge Guo \\
Department of Mathematical Sciences\\
New Jersey Institute of Technology \\
Newark, NJ 07102-1982 \\
\and
Joseph P. Romano \\
Departments of Statistics and Economics\\
 Stanford University\\
Stanford, CA 94305-4065\\
}

\date{}

\maketitle

\begin{abstract}
In this paper, we consider the problem of simultaneously testing many two-sided hypotheses
when rejections of null hypotheses are accompanied by claims of the direction of the alternative. The fundamental goal is  to construct methods that control  the mixed directional familywise error rate (mdFWER), which is the probability of making any type 1 or type 3 (directional) error.
In particular,  attention is focused on  cases where the hypotheses are ordered
as $H_1 , \ldots, H_n$, so that $H_{i+1}$ is tested only if $H_1 , \ldots, H_i$ have all been
previously rejected.
In this situation,   one can control the usual familywise error rate under arbitrary dependence by the basic procedure which tests each hypothesis at level $\alpha$, and no other multiplicity adjustment is needed. However, we show that this is far too liberal
if one also accounts for directional errors.  But, by imposing certain dependence assumptions
on the test statistics,  one can retain the basic procedure. Through a simulation study and a clinical trial example, we numerically illustrate good performance of the proposed procedures compared to the existing mdFWER controlling procedures. The proposed procedures are also implemented in the R-package FixSeqMTP.

\end{abstract}

\vspace{9pt}
\par
\vspace{6pt}

\fontsize{10.95}{14pt plus.8pt minus .6pt}\selectfont

\section{Introduction}

Directional errors or type 3 errors occur in testing situations with two-sided alternatives
when  rejections are accompanied by  additional directional claims.
For example, when testing a null hypothesis $\theta = 0$ against $\theta \ne 0$,
rejection of the null hypothesis is often augmented with the decision of whether $\theta > 0$ or
$\theta < 0$.   In the case of testing a single hypothesis, type 3 error is generally  controlled at level $\alpha$ when type 1 error is controlled at level $\alpha$ (and sometimes type 3 error is controlled at level $\alpha/2$).  However, in the case of simultaneously testing multiple hypotheses, it is often not known whether additional directional decisions can be made without losing control
of the mixed directional familywise error rate (mdFWER), the probability of at least one type 1 or type 3 error. Some methods have been developed in the literature by augmenting additional directional decisions to the existing $p$-value based stepwise procedures. Shaffer (1980) showed that Holm's procedure (Holm, 1979), augmented with decisions on direction based on the values of test statistics, can strongly control mdFWER under the assumption that the test statistics are independent and under specified conditions on the marginal distributions of the test statistics, but she also showed
that counterexamples exist even with two hypotheses.  Finner (1994) and Liu (1997) independently proved the same result for the Hochberg procedure (Hochberg, 1988). Finner (1999) generalized the result of Shaffer (1980) to a large class of stepwise or closed multiple test procedures under the same assumptions.  Some recent results have been obtained in Guo and Romano (2015).

Several situations occur in practice where hypotheses are ordered in advance, based on relative importance by some prior knowledge (for example in dose-response study, hypotheses of higher dose vs. a placebo are tested before those of lower dose vs. placebo), or  there  exists a  natural hierarchy in tested hypotheses (for example in a clinical trial, secondary endpoints are tested only when the associated primary endpoints are significant), and so on. In such fixed sequence multiple testing situations, it is also desired to make further directional decisions once significant differences are observed. For example, in dose response studies, once the hypothesis of no difference between a dose and placebo is rejected, it is of interest to decide whether the  new treatment  dose is more or less effective than the placebo. In such cases,  the possibility of making type 3 errors must be taken into account.

For control of the usual familywise error rate (FWER) (which does not account for the possibility of additional type 3 errors),  the conventional \emph{fixed sequence} multiple testing procedure that strongly controls the FWER under arbitrary dependence, is known to be a powerful procedure in testing situations with pre-ordered hypotheses (Maurer et al., 1995).  For reviews on recent relevant developments of fixed sequence multiple testing procedures for testing strictly pre-ordered hypotheses and gatekeeping strategies for testing partially pre-ordered hypotheses, see Dmitrienko, Tamhane and Bretz (2009) and Dmitrienko, Agostino and Huque (2013). Indeed, suppose null hypotheses $H_1 , \ldots , H_n$ are pre-ordered, so that $H_{i+1}$ is tested only if $H_1 , \ldots , H_i$ have all been rejected. The probability mechanism generating the data is $P$ and $H_i$ asserts that $P \in \omega_i$, some family of data generating distributions.   In such case, it is easy to see that each $H_i$ can be tested at level $\alpha$ in order to control the FWER at level $\alpha$, so that no adjustment for multiplicity is required.  The argument is simple and goes as follows.   Fix any given $P$ such that at least one $H_i$ is true (or otherwise the FWER is 0 anyway). If $H_1$ is true, i.e. $P \in \omega_1$, then a  type 1 error occurs if and only if $H_1$ is rejected, and so the FWER is just the probability $H_1$ is rejected, which is assumed controlled at level $\alpha$ when testing $H_1$.  If $ H_1$ is false, just let $f$ be the smallest index corresponding to a true null hypothesis, i.e. $H_f$ is true but $H_1 , \ldots  , H_{f-1}$ are all false.  In this case, a type 1 error occurs if and only if $H_f$ is rejected, which is assumed to be controlled at level $\alpha$.

In fact, in situations where ordering is not specified, the above result suggests it may be worthwhile to think about hypotheses in order of importance so that potentially false hypotheses are more easily detected.  Indeed, as is well-known, when the number $n$ of tested hypotheses is large, control of the FWER is often so stringent that often no rejections can be detected, largely due to the multiplicity of tests and the need to find significance at very low levels (as required, for example, in the Bonferroni method with $n$ large).   On the other hand, under a specified ordering, each test is  carried out at the same conventional level.

To our knowledge, no one explores the possibility of making additional directional decisions for such fixed sequence procedures. In this paper, we introduce such fixed sequence procedures augmented with additional directional decisions and discuss its mdFWER control under independence and some dependence. For such directional procedures, its simple fixed sequence structure of the tested hypotheses makes the notoriously challenging problem of controlling the mdFWER under dependence a little easier to handle than stepwise procedures.

 Throughout this work, we consider the problem of  testing $n$ two-sided hypotheses $H_1 , \ldots , H_n$
specified as follows:
\begin{equation}\label{eq:1}
H_i: \theta_i = 0 \quad \text{vs.} \quad H_{i}^{'}: \theta_i \neq 0, \quad i = 1, \ldots, n.
\end{equation}
We assume the hypotheses are ordered in advance, either using some prior knowledge about the importance of the hypotheses or by some other specified criteria, so that $H_1$ is tested first
and $H_i$ is only tested if $H_1 , \ldots , H_{i-1}$ are all rejected.
We also assume that, for each $i$, a  test statistic $T_i$ and $p$-value $P_i$ are available to test $H_i$
(as a single test).
For a rejected hypothesis $H_i$, we decide on the sign of the parameter $\theta_i$ by the sign of the corresponding test statistic $T_i$, i.e., we conclude $\theta_i > 0$ if $T_i > 0$ and vice versa. The errors that might occur while testing these hypotheses are type 1 and type 3 errors. A \emph{type 1 error} occurs when a true $H_i$ is falsely rejected. A \emph{type 3 error} occurs when a false $H_i$ is correctly rejected but  the claimed sign of the parameter $\theta_i$ is wrong.
Then, the mdFWER is the probability of making at least a type 1 or type 3 error,
and it is desired that this error rate is no bigger than $\alpha$ for all possible data generating
distributions in the model.

We make a few standard assumptions about the test statistics. Let $T_i \sim F_{\theta_i}(\cdot)$ for some continuous cumulative distribution function  $F_{\theta_i} ( \cdot ) $ having parameter $\theta_i$.
In general, most of our results also apply through the same arguments when the family of distributions of $T_i$  depends on $i$, though
for simplicity of notation,  the notation is suppressed. We assume that $F_0$ is symmetric about $0$ and $F_{\theta_i}$ is stochastically increasing in $\theta_i$.   Various dependence assumptions  between the test statistics will be used throughout
the paper.  (Some of the results can generalize outside this parametric framework.
Of course, for many problems, approximations are used to construct marginal tests and the
approximate distributions of the $T_i$ are often normal, in which case our exact finite sample results will hold approximately as well.)
Let $c_1 = F_0^{-1} (  \alpha / 2 )$ and $c_2 = F_0^{-1} ( 1- \alpha /2 )$, so that a marginal
level $\alpha $ test of $H_i$ rejects if $T_i < c_1$ or $T_i > c_2$.
For testing $H_i$ vs. $H_{i}^{'}$,   rejections are based on large values of $|T_i|$ and the corresponding two-sided $p$-value is defined by
\begin{equation}\label{eq:3}
P_i = 2\min \{ F_{0}(T_i), 1-F_{0}(T_i) \}, \quad i = 1,\ldots, n.
\end{equation}
We assume that the  $p$-value $P_i$ is distributed as U(0,1) when $\theta_i = 0$.

The rest of the paper is organized as follows.
In Section \ref{section:arbitrary},
we consider the problem of  mdFWER control under no dependence assumptions
on the test statistics.
Unlike control of the usual FWER where each test can be constructed at level $\alpha$, it
is seen that $H_i$ can only be tested at a much smaller level $\alpha / 2^{i-1}$.
This rapid decrease in the critical values used motivates the study of the problem
under various dependence assumptions.
In Section 3 we introduce a directional fixed sequence procedure and prove that this procedure controls the mdFWER under independence. In Sections 4 and 5 we further discuss its mdFWER control under positive dependence.  In Section 6 we numerically evaluate the performances of the proposed procedure through a simulation study. In Section 7  we illustrate an application of the proposed procedures through a clinical trial example. Section 8 makes some concluding remarks and all proofs are deferred to Section 9.

\section{The mdFWER Control Under Arbitrary Dependence}\label{section:arbitrary}

A general fixed sequence procedure based on marginal $p$-values must specify
the critical level $\alpha_i$ that is used for testing $H_i$, in order for the resulting
procedure to control the mdFWER at  level $\alpha$.    When controlling the FWER without
regard to type 3 errors, each $\alpha_i$ can be as large as $\alpha$.  However,
Theorem \ref{theorem:nodep}  below shows that  by  using the critical constant $\alpha_i = \alpha/2^{i-1}$, the mdFWER is controlled at level $\alpha$.  Moreover, we show
that these critical constants are unimprovable.
Formally, the optimal procedure is defined as follows.

\begin{procedure}[\textbf{Directional fixed sequence procedure under arbitrary dependence}] \hfill
\begin{itemize}
\item Step 1: If $P_1 \le \alpha$ then reject $H_1$ and continue to test $H_2$ after making directional decision on $\theta_1$: conclude $\theta_1 > 0$ if $T_1 > 0$ or  $\theta_1 <0$ if $T_1 < 0$. Otherwise, accept all the hypotheses and stop.
\item Step $i$:  If $P_i \le \alpha/2^{i-1}$ then reject $H_i$ and continue to test $H_{i+1}$ after making directional decision on $\theta_i$: conclude $\theta_i > 0$ if $T_i > 0$ or $\theta_i <0$ if $T_i < 0$. Otherwise, accept the remaining  hypotheses $H_i, \ldots, H_n$.
\end{itemize}\label{procedure:nodep}
\end{procedure}

In the following, we discuss the mdFWER control of Procedure \ref{procedure:nodep} under arbitrary dependence of the $p$-values. When testing a single hypothesis, the mdFWER of Procedure \ref{procedure:nodep} reduces to the type 1 or type 3 error rate depending on whether $\theta = 0$ or $\theta \neq 0$, and Procedure \ref{procedure:nodep} reduces to the usual $p$-value based method along with the directional decision for the two-sided test. The following lemma covers this case.

\begin{lemma}\label{lem1}
Consider testing the single hypothesis $H: \theta = 0$ against $H^{'}: \theta \neq 0$ at level $\alpha$, using the usual $p$-value based method along with a directional decision. If $H$ is a false null hypothesis, then the type 3 error rate is bounded above by $\alpha/2$.
\end{lemma}

Generally, when simultaneously testing $n$ hypotheses, by using Lemma \ref{lem1} and mathematical induction, we have the following result holds.

\begin{theorem}\label{theorem:nodep}
For Procedure \ref{procedure:nodep} defined as above, the following conclusions hold.
\begin{itemize}
    \item[(i)]This procedure
strongly controls the mdFWER at level $\alpha$ under arbitrary
dependence of the $p$-values.
    \item[(ii)] One cannot increase even one of the critical constants
    $\alpha_i = \alpha/2^{i-1}, i = 1, \ldots, n,$ while keeping the remaining fixed without losing control of the mdFWER.
\end{itemize}
\end{theorem}

In fact, the proof shows that no strong parametric assumptions are required.  However,
the rapid decrease in critical values $\alpha / 2^{i-1}$ makes rejection of additional hypotheses
difficult.  Thus, it is of interest to explore how dependence assumptions can be used
to increase these critical constants while maintaining control of the mdFWER.
The assumptions and methods will be described in the remaining sections.

\begin{remark} \rm
Instead of Procedure \ref{procedure:nodep}, let us consider the conventional fixed sequence procedure with the same critical constant $\alpha$ augmented with additional directional decisions, which is defined in Section 3 as Procedure 2. By using Bonferroni inequality and Lemma 1, we can prove that the mdFWER of this procedure is bounded above by $\frac{n+1}{2} \alpha$. Thus, the modified version of the procedure, which has the same critical constant $\frac{2\alpha}{n+1}$, strongly controls the mdFWER at level $\alpha$ under arbitrary dependence of $p$-values. However, it is unclear if such critical constant can be further improved without losing the control of the mdFWER.
\end{remark}

\section{The mdFWER Control Under Independence}\label{section:independence}

We further make the following assumptions on the distribution of the test statistics.

\begin{assumption}[\textbf{Independence}]\label{assumption:indep}
The test statistics, $T_{1}, \ldots, T_{n}$, are mutually independent. \end{assumption}

Of course, it follows that the $p$-values $P_{1}, \ldots, P_{n}$ are mutually independent as well.

As will be seen, it will be necessary to make further assumptions on the family of distributions
for each marginal test statistic.

\begin{definition}[\textbf{Monotone Likelihood Ratio (MLR)}]
A  family of probability density functions  $f_{\delta}(\cdot)$ is said to have monotone likelihood ratio property if, for any two values of the parameter $\delta$, $\delta_2 > \delta_1$ and any two points $x_2 > x_1$,
\begin{equation}\label{eq:4}
\frac{f_{\delta_2}(x_2)}{f_{\delta_1}(x_2)} \ge \frac{f_{\delta_2}(x_1)}{f_{\delta_1}(x_1)},
\end{equation}
or equivalently,
\begin{equation}\label{eq:5}
\frac{f_{\delta_1}(x_1)}{f_{\delta_1}(x_2)} \ge \frac{f_{\delta_2}(x_1)}{f_{\delta_2}(x_2)}.
\end{equation}
\end{definition}

Definition 1 means that, for fixed $x_1 < x_2$, the ratio $\frac{f_{\delta}(x_1)}{f_{\delta}(x_2)}$ is non-increasing in $\delta$. Two direct implications of Definition 1 in terms of the cdf $F_{\delta}(\cdot)$ are
\begin{equation}\label{eq:6}
\frac{F_{\delta_1}(x_2)}{F_{\delta_1}(x_1)} \le \frac{F_{\delta_2}(x_2)}{F_{\delta_2}(x_1)},
\end{equation}
and
\begin{equation}\label{eq:7}
\frac{1-F_{\delta_1}(x_2)}{1-F_{\delta_1}(x_1)} \le \frac{1-F_{\delta_2}(x_2)}{1-F_{\delta_2}(x_1)}.
\end{equation}

\begin{assumption}[MLR Assumption]\label{assumption:mlr}
The family of   marginal distributions of the $T_i$  has   monotone likelihood ratio.
\end{assumption}

Based on the conventional fixed sequence multiple testing procedure, we define a directional fixed sequence procedure as follows, which is the conventional fixed sequence procedure augmented with directional decisions.  In other words, any hypothesis is tested at level $\alpha$, and as will be seen under the specified conditions,  no
reduction in critical values is necessary in order to achieve mdFWER control.

\begin{procedure}[\textbf{Directional fixed sequence procedure}] \hfill
\begin{itemize}
\item Step 1: If $P_1 \le \alpha$, then reject $H_1$ and continue to test $H_2$ after making a directional decision on $\theta_1$: conclude $\theta_1 > 0$ if $T_1 > 0$ or  $\theta_1 <0$ if $T_1 <  0$. Otherwise, accept all the hypotheses and stop.
\item Step $i$:  If $P_i \le \alpha$ , then reject $H_i$ and continue to test $H_{i+1}$ after making a directional decision on $\theta_i$: conclude $\theta_i > 0$ if $T_i > 0$ or $\theta_i <0$ if $T_i <  0$. Otherwise, accept the remaining  hypotheses, $H_i, \ldots, H_n$.
\end{itemize}\label{procedure:indep}
\end{procedure}

For Procedure \ref{procedure:indep}, in the case of $n=2$,  we  derive a simple expression for the mdFWER in Lemma \ref{lem2} below and prove its mdFWER control in Lemma \ref{lem3} by using such simple expression.

\begin{lemma}\label{lem2}
Consider testing two hypotheses $ H_1:\theta_1 = 0$ and $H_2: \theta_2 = 0$, against both sided alternatives, using Procedure \ref{procedure:indep} at level $\alpha$. Let $c_1 = F_{0}^{-1}(\alpha/2 )$ and $c_2 = F_{0}^{-1}(1- \alpha/2)$. When $\theta_2 = 0$, the following result holds.
\begin{equation}\label{equation:old8}
mdFWER = \left\{
  \begin{array}{l l}
    \alpha + F_{\theta_1}(c_1) - F_{\theta_1}(c_2) + F_{(\theta_1,0)}(c_2, c_2) - F_{(\theta_1,0)}(c_2, c_1) & \text{if $\theta_1 >  0$}\\
    1 +  F_{\theta_1}(c_1) - F_{\theta_1}(c_2) + F_{(\theta_1, 0)}(c_1, c_1) - F_{(\theta_1, 0)}(c_1, c_2) & \text{if $\theta_1 <  0$.}
  \end{array} \right.
\end{equation}
\noindent   In the above, $F_{\theta_1 , \theta_2 } ( \cdot , \cdot )$ refers to the joint c.d.f. of $(T_1 , T_2 )$.  Then, under Assumption \ref{assumption:indep} (independence) , (\ref{equation:old8}) can be simplified as
\begin{equation}\label{equation:old9}
mdFWER = \left\{
  \begin{array}{l l}
    \alpha + F_{\theta_1}(c_1) - \alpha F_{\theta_1}(c_2) & \text{if $\theta_1 >  0$}\\
    1 + \alpha F_{\theta_1}(c_1) - F_{\theta_1}(c_2) & \text{if $\theta_1 <  0$.}
  \end{array} \right.
\end{equation}
\end{lemma}

\begin{lemma}\label{lem3}
Under Assumption \ref{assumption:indep} (independence) and Assumption \ref{assumption:mlr} (MLR), Procedure \ref{procedure:indep} strongly controls the mdFWER when $n=2$.
\end{lemma}

Generally, for testing any $n$ hypotheses, by using mathematical induction and Lemma \ref{lem3}, we also prove the mdFWER control of Procedure \ref{procedure:indep} under the same assumptions as in the case of $n=2$.

\begin{theorem}\label{theorem:indep}
Under Assumption 1  (independence)  and Assumption 2 (MLR), Procedure \ref{procedure:indep} strongly controls the mdFWER at level $\alpha$.
\end{theorem}

Many families of distributions have the MLR property: normal, uniform, logistic, Laplace, Student's t, generalized extreme value, exponential  familes of distributions, etc.
However, it is also important to know whether or not the above results fail without the MLR assumption.
A natural family of distributions to consider without the MLR property is the Cauchy family; indeed,
Shaffer (1980) used this family to obtain a counterexample for the directional Holm procedure while testing $p$-value ordered hypotheses. We now show that Procedure \ref{procedure:indep} fails to  control  the mdFWER for this family of distributions with corresponding cdf $F_{\theta}(x) = 0.5 + \frac{1}{\pi} \arctan(x - \theta)$, even under independence.

Lemma \ref{lem2} can  be used to verify the calculation for the case of $n=2$ with $\theta_1 > 0$ and $\theta_2 = 0$; specifically, see (\ref{equation:old9}). Indeed,
we just need to show  show
\begin{equation}\label{equation:violate}
F_{\theta_1}(-c) = F_0(-c - \theta_1) > \alpha F_{\theta_1}(c) = \alpha F_0(c - \theta_1)~,
\end{equation}
 where $c$ is the $1 - \alpha/2$ quantile
of the standard Cauchy distribution, given by $\tan [\pi (1-\alpha)/2]$. Take $\alpha = 0.05$, so $c = 12.7062$. Then, the above inequality (\ref{equation:violate}) is violated for example by $\theta_1 = 100$. The left side is approximately $F(-112.7) \approx 0.002824$ while the right side is
$$0.05 \times F(-87.3) = 0.05 \times 0.0036 = 0.00018.$$

\section{Extension to Positive Dependence}\label{section:positive}

Clearly, the assumption of independence is of limited
utility in multiple testing, as many tests are usually carried out on the same data set.
Thus, it is important to generalize the results of the previous section to cover some more general cases.  As is typical in the multiple testing literature (Benjamini and Yekutieli, 2001; Sarkar, 2002; Sarkar and Guo, 2010, etc), assumptions of positive regression dependence will be used.

Before defining the assumptions, for convenience, we introduce several notations below. Among the prior-ordered hypotheses $H_1, \ldots, H_n$, let $i_0$
denote the index of the first true null hypothesis, $n_1$ denote the number of all false nulls,
and $T_{i_1}, \ldots, T_{i_{n_1}}$ denote the corresponding false null test statistics. Specifically, if all $H_i$'s are false, let $i_0 = n+1$.

\begin{assumption}\label{assumption:prd}
The false null test statistics along with parameters, $\theta_{i_1} T_{i_1}, \ldots, \theta_{i_{n_1}} T_{i_{n_1}}$, are positively regression dependent
in the sense of
\begin {eqnarray}
E \left \{ \phi(\theta_{i_1} T_{i_1}, \ldots, \theta_{i_{n_1}} T_{i_{n_1}}) ~|~ \theta_{i_k}T_{i_k} \ge u \right \} \uparrow u,
\end {eqnarray}
for each $\theta_{i_k}T_{i_k}$ and any (coordinatewise) non-decreasing function
$\phi$.
\end{assumption}

\begin{assumption}\label{assumption:indep2}
The first true null statistic, $T_{i_0}$, is independent of all false null statistics
$T_{i_k}, k=1, \ldots, n_1$ with $i_k  <  i_0$.
\end{assumption}

\begin{theorem}\label{theorem:positive}
Under Assumptions \ref{assumption:mlr} - \ref{assumption:indep2}, Procedure \ref{procedure:indep} strongly controls the mdFWER at level $\alpha$.
\end{theorem}

\begin{corollary}\label{corollary:1}
When all tested hypotheses are false, Procedure \ref{procedure:indep} strongly controls the mdFWER at level $\alpha$ under Assumptions \ref{assumption:mlr} - \ref{assumption:prd}.
\end{corollary}

\begin{remark} \rm
In Theorem \ref{theorem:positive}, we note that specifically, when all of the tested hypotheses are false, Assumption \ref{assumption:indep2} is automatically satisfied. Generally, consider the case of any combination of true and false null hypotheses where Assumption \ref{assumption:indep2} is not imposed.
Without loss of generality, suppose $\theta_i > 0, i=1, \ldots, n-1$ and $\theta_n = 0$, that is, the first $n-1$ hypotheses are false and the last one is true. Under Assumptions \ref{assumption:mlr}-\ref{assumption:prd}, if $T_n$ (or $-T_n$) and $T_1, \ldots, T_{n-1}$ are positively regression dependent, then the mdFWER of Procedure \ref{procedure:indep} when testing $H_1, \ldots, H_n$ is, for any $n$,  bounded above by
\begin{eqnarray}
&& \text{Pr( make at least one type 3 error when testing }H_1, \ldots, H_{n-1} \text{ or } T_n \notin (c_1, c_2)) \nonumber \\
& \le & \lim_{\theta_n \rightarrow 0+}\text{Pr( make at least one type 3 error when testing }H_1, \ldots, H_n) \nonumber \\
&& \quad + ~\lim_{\theta_n \rightarrow 0+}Pr(T_n \ge c_2) \nonumber \\
& \le & \alpha + \alpha/2 = 3\alpha/2. \nonumber
\end{eqnarray}
The first inequality follows from the fact that when $\theta_n \rightarrow 0+$, $H_n$ can be interpreted as a false null hypothesis with $\theta_n > 0$, and thus one type 3 error is made if $H_n$ is rejected and $T_n \le c_1$. The second inequality follows from Corollary \ref{corollary:1}
and Lemma \ref{lem1}.

Based on the above inequality, a modified version of Procedure \ref{procedure:indep}, the directional fixed sequence procedure with the critical constant $2\alpha/3$, strongly controls the mdFWER at level $\alpha$ under Assumptions \ref{assumption:mlr}-\ref{assumption:prd} and the above additional assumption.
\end{remark}

\begin{remark} \rm
In the above remark, further, if we do not make any assumption regarding dependence between the true null statistic $T_n$ and the false null statistics $T_1, \ldots, T_{n-1}$. Then, by Theorem \ref{theorem:positive}, the mdFWER of Procedure \ref{procedure:indep} when testing $H_1, \ldots, H_n$ is bounded above by
\begin{eqnarray}
&& \text{Pr( make at least one type 3 error when testing }H_1, \ldots, H_{n-1}) \nonumber \\
&& \quad + ~\text{Pr( make type 1 error when testing } H_n) \nonumber \\
& \le & \alpha + \alpha = 2\alpha. \nonumber
\end{eqnarray}
Therefore, an alternative modified version of Procedure \ref{procedure:indep}, the directional fixed sequence procedure with the critical constant $\alpha/2$, strongly controls the mdFWER at level $\alpha$ only under Assumptions \ref{assumption:mlr}-\ref{assumption:prd}.
\end{remark}

\section{Further Extensions to Positive Dependence}\label{section:positive}

We now develop alternative results
 to show that Procedure \ref{procedure:indep} can control mdFWER even under certain dependence between the false null and true null statistics. We relax the assumption of independence that the false null statistics are independent of the first true null statistic, and consider a slightly strong version of the conventional positive regression dependence on subset of true null statistics (PRDS) (Benjamini and Yekutieli, 2001), which is given below.

\begin{assumption}\label{assumption:prds2}
The false null test statistics, $T_1, \ldots, T_{i_0-1}$
and the first true null statistic $T_{i_0}$, are positive regression dependent
in the sense of
\begin {eqnarray}
E \left \{ \phi(T_1, \ldots, T_{i_0-1}) ~|~ T_{i_0} \ge u, T_1, \ldots, T_j \right \} \uparrow u,
\end {eqnarray}
for any given $j=1, \ldots, i_0 -1$, any given values of $T_1, \ldots, T_j$ and any (coordinatewise) non-decreasing function
$\phi$.
\end{assumption}

We firstly consider the case of $n=2$, that is, while testing two hypotheses, and show control of the mdFWER of Procedure \ref{procedure:indep} when the test statistics are positively regression dependent in the sense of  Assumption \ref{assumption:prds2}.

\begin{proposition}\label{proposition:1}
Under Assumptions \ref{assumption:mlr} and \ref{assumption:prds2}, the mdFWER of Procedure \ref{procedure:indep} is strongly controlled at level $\alpha$ when $n=2$.
\end{proposition}

Specifically, in the case of bivariate normal distribution, Assumption \ref{assumption:mlr} is satisfied and two test statistics $T_1$ and $T_2$ are always positively or negatively regression dependent. As in the proof of Proposition \ref{proposition:1}, to show the mdFWER control of  Procedure \ref{procedure:indep}, we only need to consider the case of $\theta_1 \neq 0$ and $\theta_2 = 0$. Thus, if $T_1$ and $T_2$ are negatively regression dependent, we can choose $-T_2$ as the statistic for testing $H_2$ and Assumption \ref{assumption:prds2} is still satisfied. By Proposition \ref{proposition:1}, we have the following corollary holds.

\begin{corollary}
Under the case of bivariate normal distribution, the mdFWER of Procedure \ref{procedure:indep} is strongly controlled at level $\alpha$ when $n=2$.
\end{corollary}

We now consider the case of three hypotheses. The general case will ultimately be considered, but is instructive to discuss the case separately due to the added multivariate MLR condition, which is described as follows.

Let $f(x|T_1)$ and $g(x|T_1)$ denote the probability density functions of $T_2$ and $T_3$ conditional on $T_1$, respectively.

\begin{assumption}[\textbf{Bivariate  Monotone Likelihood Ratio (BMLR)}]\label{assumption:bmlr}
For any given value of $T_1$, $f(x|T_1)$ and $g(x|T_1)$ have the monotone likelihood ratio (MLR) property in $x$, i.e., for any $x_2 > x_1$, we have
\begin{eqnarray}\label{eq:22}
\frac{f(x_2|T_1)}{g(x_2|T_1)} \ge \frac{f(x_1|T_1)}{g(x_1|T_1)}.
\end{eqnarray}
\end{assumption}

\begin{proposition}\label{proposition:2}
Under Assumptions \ref{assumption:mlr}, \ref{assumption:prd}, \ref{assumption:prds2}, and \ref{assumption:bmlr}, the mdFWER of Procedure \ref{procedure:indep} is strongly controlled at level $\alpha$ when $n=3$.
\end{proposition}

\begin{remark}\rm
In the case of three hypotheses, suppose that the test statistics $T_i, i=1, \ldots, 3$ are trivariate normally distributed with the mean $\theta_i$. Without loss of generality, assume $\theta_i > 0, i=1, 2$ and $\theta_3 =0$, that is, $H_1$ and $H_2$ are false and $H_3$ is true. Let $\Sigma = (\sigma_{ij}), i, j=1, \ldots, 3$, denote the variance-covariance matrix of $T_i$'s. It is easy to see that Assumption \ref{assumption:mlr}  is always satisfied. Also, when $\sigma_{ij} \ge 0$ for $i \neq j$, Assumption \ref{assumption:prd} and Assumption \ref{assumption:prds2} are satisfied. Finally, when $\sigma_{22} = \sigma_{33}$ and $\sigma_{12} = \sigma_{13}$, Assumption \ref{assumption:bmlr}  is satisfied.
\end{remark}

\bigskip

Finally, We consider the general  case of $n$ hypotheses.   Now we must consider
the multivariate monotone likelihood ratio property, described as follows.
For any given $j=1, \ldots, i_0-1$, let $f(x|T_1, \ldots, T_{j-1})$ and $g(x|T_1, \ldots, T_{j-1})$ denote the probability density functions of $T_j$ and $T_{i_0}$ conditional on $T_1, \ldots, T_{j-1}$, respectively.

\begin{assumption}[\textbf{Multivariate  Monotone Likelihood Ratio (MMLR)}]\label{assumption:mmlr}
For any given values of \\ $T_1, \ldots, T_{j-1}$,  $f(x|T_1, \ldots, T_{j-1})$ and $g(x|T_1, \ldots, T_{j-1})$ have the monotone likelihood ratio (MLR) property in $x$, i.e., for any $x_2 > x_1$, we have
\begin{eqnarray}\label{eq:32}
\frac{f(x_2|T_1, \ldots, T_{j-1})}{g(x_2|T_1, \ldots, T_{j-1})} \ge \frac{f(x_1|T_1, \ldots, T_{j-1})}{g(x_1|T_1, \ldots, T_{j-1})}.
\end{eqnarray}
\end{assumption}

\begin{theorem}\label{theorem:n}
Under Assumptions \ref{assumption:mlr}, \ref{assumption:prd}, \ref{assumption:prds2}, and \ref{assumption:mmlr}, the mdFWER of Procedure \ref{procedure:indep} is strongly controlled at level $\alpha$.
\end{theorem}

\section{A Simulation Study}

We conduct a simulation study to illustrate the performance of the proposed directional fixed sequence procedures under arbitrary dependence (Procedure 1) and independence (Procedure 2) in terms of mdFWER control and average power and compare them with the directional Bonferroni procedure, directional Holm procedure and directional Hochberg procedure. We study two simulation settings for evaluating the effects of proportion of false nulls and dependence on the performance of these procedures, respectively. We generate $n$-dimensional normal random vectors $(T_1, \ldots, T_n)$ where the components follow normal distribution $N(\theta_i, 1)$ with common pairwise correlation $\rho$.
Consider simultaneously testing $n$ two-sided hypotheses using $T_i$ along with making directional decisions on $\theta_i$ based on the sign of $T_i$:
\begin{equation}\label{eq:1}
H_i: \theta_i = 0 \quad \text{vs.} \quad H_{i}^{'}: \theta_i \neq 0, \quad i = 1, \ldots, n.
\end{equation}
For this simulation, we set the first $n_1$ of the $n$ hypotheses $H_i$ to be false null and the rest to be true null. The true null test statistics are generated from $N(0, 1)$ and the false null test statistics are generated from $N(\theta_i, 1)$ with $\theta_i \neq 0$. The simulation results are obtained under the significance level $\alpha=0.05$ and based on 10,000 replicates. The ``power'' of a procedure at a replication is defined as the proportion of non-null $\theta_i$ to be rejected along with correct directional decisions on $\theta_i$ to be made among all non-null $\theta_i$ out of $n$ hypotheses. The ``average power'' is the average of the power for the 10,000 replications. The mdFWER is estimated as the proportion of replications where at least one true null hypothesis is falsely rejected or at least one false null hypothesis is correctly rejected but a wrong directional decision is made regarding
the corresponding $\theta_i$.

\subsection{Simulation Setting 1}

In this setting, we set the number of tested hypotheses $n = 20$,  the common correlation $\rho = 0$ (independence) or $\rho = 0.5$ (positive correlation), and the proportion of false null hypotheses $\pi_1$ to be between $0.05$ and $1.0$. For the values of $\theta_i$, we set $\theta_i = 3$ for false null and $\theta_i = 0$ for true null.

\begin{figure}[h!]
  \includegraphics[width=\linewidth]{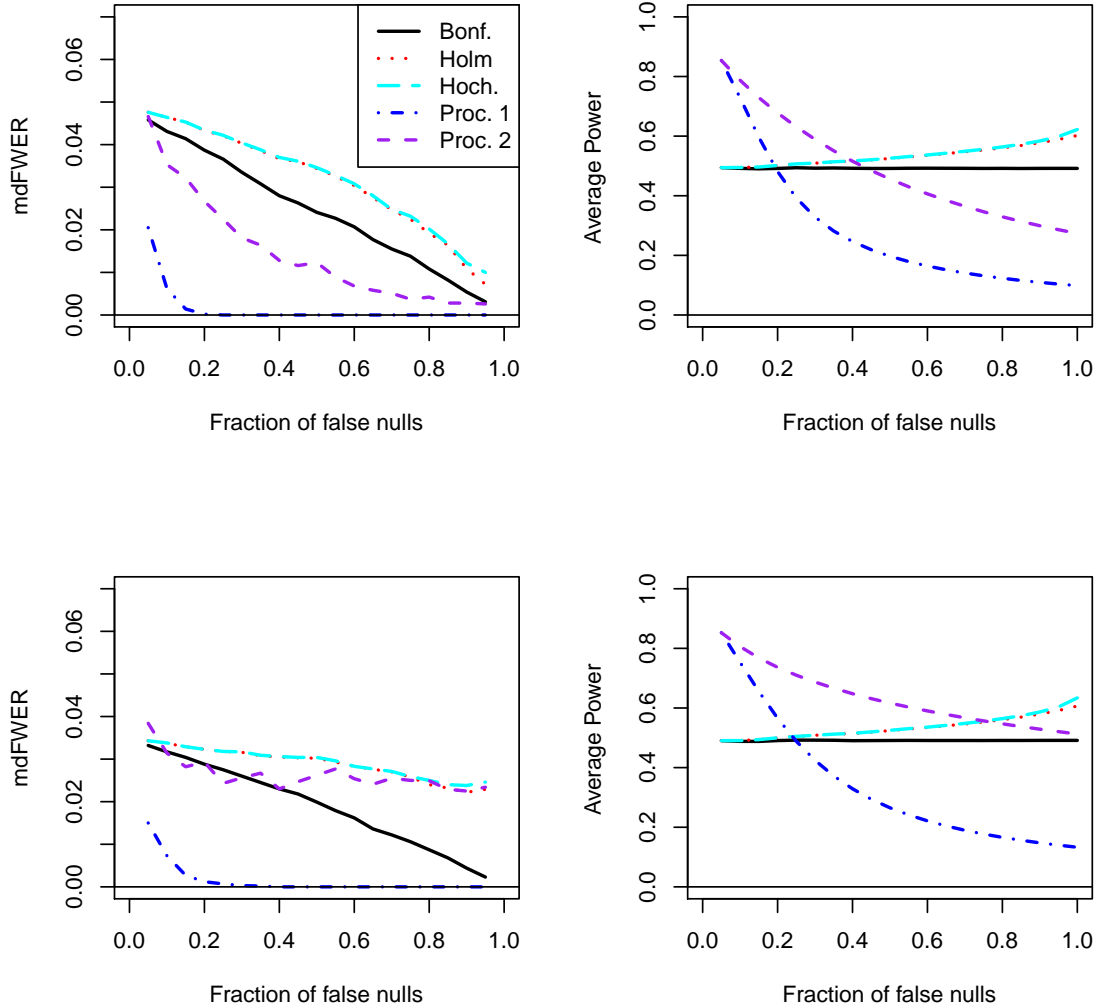}
  \caption{Estimated mdFWER and average powers of our suggested Procedure 1 (Proc. 1) and Procedure 2 (Proc. 2) along with existing directional Bonferroni procedure (Bonf.), directional Holm procedure (Holm), and directional Hochberg procedure (Hoch.) for $n=20$ hypotheses with the fraction of false nulls $\pi_1$ from $0.05$ to $1.0$ and common correlation $\rho = 0$ (upper panel) or $\rho = 0.5$ (bottom panel).}\label{fig:Setting1}
\end{figure}

Figure \ref{fig:Setting1} shows the plots of mdFWER and average power of all five directional procedures plotted against $\pi_1$, the fraction of false null hypotheses. As it is evident, all the five procedures control mdFWER at level 0.05 and Procedure 1 has the lowest mdFWER. When the test statistics are independent ($\rho = 0$), the mdFWER of Procedure 2 is also lower than those of the existing procedures, whereas when the test statistics are positively correlated ($\rho = 0.5$), the mdFWER of Procedure 2 is generally higher than that of the directional Bonferorni procedure but lower than those of the directional Holm and directional Hochberg procedures except for very high fractions of false nulls.

When the fraction of false nulls is low or moderate ($\pi_1 \le 0.4$), as is usually expected in practical applications, Procedure 2 has the highest power followed by Procedure 1, both when the test statistics are independent or positively correlated. However,
when the fraction of false nulls is high, even Procedure 2 loses its edge over the existing procedures. We also observe from Figure \ref{fig:Setting1} that the proposed procedures and the existing procedures have different power performances with increasing proportion of false nulls.
The average powers of Procedures 1 and 2 are decreasing in terms of the proportion of false nulls, whereas the average powers of the existing Procedures are slightly increasing in the proportion of false nulls.

\begin{figure}[h!]
  \includegraphics[width=\linewidth]{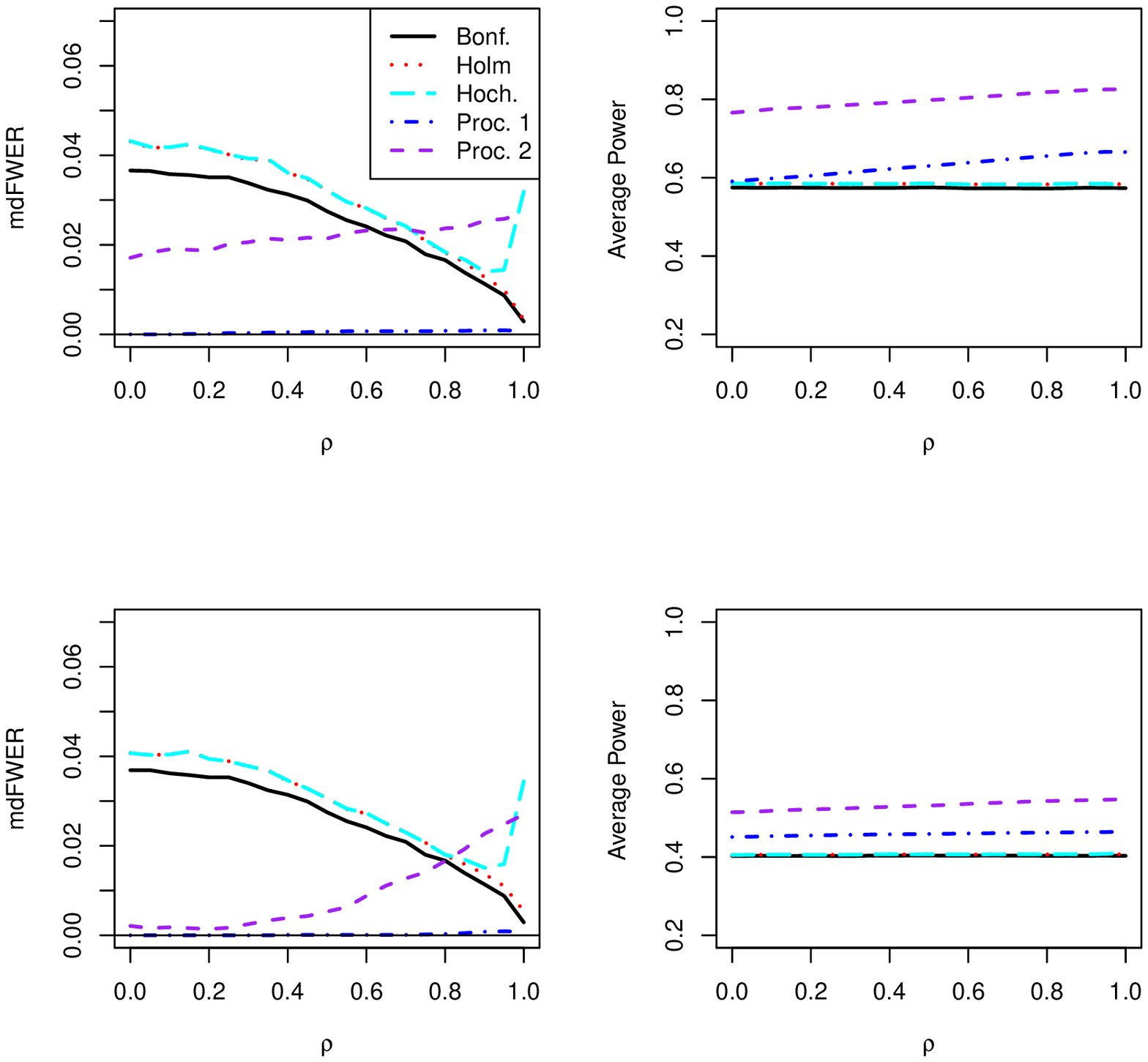}
  \caption{Estimated mdFWER and average powers of our suggested Procedure 1 (Proc. 1) and Procedure 2 (Proc. 2) along with existing directional Bonferroni procedure (Bonf.), directional Holm procedure (Holm), and directional Hochberg procedure (Hoch.) for $n=20$ hypotheses with common correlation $\rho$ between $0$ and $1$ and $n_1 = 5$ non-null $\theta_i = \theta_0 r^{i-1}$ with $(\theta_0, r) = (5, 0.8)$ (upper panel) or $(\theta_0, r) = (8, 0.5)$ (bottom panel).}\label{fig:Setting2}
\end{figure}

\subsection{Simulation Setting 2}

In this setting, we set the number of tested hypotheses $n = 20$,  the number of false null hypotheses
$n_1 = 5$, and the common correlation $\rho$ to be between $0$ and $1$. For the values of non-null $\theta_i$, we set $\theta_i = \theta_0 r^{i-1}, i = 1, \ldots, n_1$, which are decreasing proportionally with the values of parameters $(\theta_0, r) = (5, 0.8)$ or
$(\theta_0, r) = (8, 0.5)$, and for the values of null $\theta_i$, we set $\theta_i = 0$.

Figure \ref{fig:Setting2} shows the plots of mdFWER and average power of all five directional procedures plotted against $\rho$, the common correlation. As seen from Figure \ref{fig:Setting2}, all the five procedures control the mdFWER at level $\alpha$ and Procedure 2 has the highest average power followed by Procedure 1 for different values of $\rho$. We also observe that our proposed procedures have different behaviors of performance with respect to common correlation compared to the existing procedures. The mdFWER and average powers of Procedures 1 and 2, and their power improvements over the existing three procedures are all increasing in terms of correlation, whereas the mdFWERs of the existing three procedures are basically decreasing in terms of correlation, except for the directional Hochberg procedure, its mdFWER becomes to be increasing when $\rho$ is very large.

\section{Clinical Trial Example}

The directional fixed sequence procedure comes in handy in dose-response studies or studies with multiple endpoints where hypotheses are ordered in advance. To illustrate our procedure we use the hypertension trial example considered in Dmitrienko et al. (2005, Page 118). This clinical trial was conducted to test the efficacy and safety of four doses of an investigational drug versus placebo. The four doses, from lowest to highest, were respectively labeled as D1, D2, D3 and D4 and the placebo was labeled P. The primary efficacy endpoint was the reduction in diastolic blood pressure (measured in mm Hg). Dose D4 was believed to be the most efficacious one (in terms of its effect on diastolic blood pressure), followed by doses D3 and D2 and dose D1 was expected to be marginally efficacious.

The original analysis had 8 two sided hypotheses, four were dose-placebo contrasts and four dose-dose contrasts. For our analysis, we use these comparisons to test the hypotheses in the order mentioned and conclude on the direction of efficacy. We apply the directional fixed-sequence procedures (Procedures 1 and 2) described in the paper and for comparison, we also include the results of the Bonferroni single-step procedure appended with directional decisions. Table \ref{table_dataanalysis} shows the results of our analysis done at level $\alpha=0.05$.

\begin{table}
\caption{Results of Directional Fixed Sequence Procedures in the hypertension trial example (P=Placebo and D1-D4 denote four doses of the investigational drug). The overall Type I error rate $\alpha = 0.05$. R: Rejected, NR: Not rejected.} \label{table_dataanalysis}
\begin{center}
\begin{tabular}{ |c|c|c|c|c|c| }
 \hline
Test & Test & Raw & {Procedure 1} & {Procedure 2} & {Bonferroni} \\
Contrast & statistic & $p$-value & Decision  & Decision  & Decision \\
 & & & (Direction) & (Direction) & (Direction) \\
\hline
 D4-P & 3.4434 & 0.0008 & R & R & R \\
 & & & (More Effective) & (More Effective) & (More Effective) \\
 D3-P & 2.5085 & 0.0135 & R & R & R \\
& & & (More Effective) & (More Effective) & (More Effective) \\
 D2-P & 2.3642 & 0.0197 & NR & R & NR \\
& & & & (More Effective) & \\
 D1-P & -0.3543 & 0.7237 & -- & NR & NR \\
 & & & & & \\
D4-D1 & 3.7651 & 0.0003 & -- & -- & R \\
 & & & & & (More Effective) \\
 D4-D2 & 1.0900 & 0.2779 & -- & -- & NR \\
 & & & & & \\
D3-D1 & 2.8340 & 0.0054 & -- & -- & R \\
 & & & & & (More Effective) \\
D3-D2 & 0.1930 & 0.8473 & -- & -- & NR \\
 & & & & & \\
 \hline
\multicolumn{3}{|c|}{Number Rejected} & {2} & {3} & {4} \\
\hline
\end{tabular}
\end{center}
\end{table}

As seen in Table \ref{table_dataanalysis} even though the Bonferroni single-step procedure appended with directional decision rejects the most number of hypotheses, two of the rejected hypotheses (D4-D1 an D3-D1) do not make sense as the hypothesis comparing D1-P is not rejected. As a result, the proposed Procedure 2 performs best rejecting the highest number of hypotheses at level 0.05 and also the analysis results make valid conclusions. We can conclude that doses 4, 3 and 2 are significantly more effective than the placebo but Dose 1 is not significantly different than the placebo.

While Procedure 1 assumes nothing about the dependence structure of the $p$-values, it is obviously then more conservative than Procedure 2. However, in the context here, where each $p$-value corresponds to a different dose of the same drug, it is reasonable to assume positive dependence of the outcomes. In such case, results based on Procedure 2 are valid and indicate that even Dose 2 is significantly beneficial as compared to the placebo.

\section{Conclusions}

In this paper, we consider the problem of simultaneously testing multiple prior-ordered hypotheses accompanied by directional decisions. The conventional fixed sequence procedure augmented with additional directional decisions are proved to control the mdFWER under independence and some dependence, whereas, this procedure is also shown to be far too liberal to control the mdFWER, if no dependence assumptions are imposed on the test statistics. Through a simulation study, we numerically show the good performances of the proposed procedures in terms of the mdFWER control and average power as compared to several existing directional procedures, directional Bonferroni, Holm, and Hochberg procedures. The proposed procedures are also implemented in the R-package FixSeqMTP.

We need to note that in the existing literature, to our knowledge, only directional Bonferroni procedure is theoretically proved to strongly control the mdFWER under dependence. It is still an open problem that the directional Holm and Hochberg procedures control the mdFWER under certain dependence. Our suggested directional fixed sequence procedure can be a powerful alternative solution to the problem of directional errors control under dependence. We hope that the approaches and techniques developed in this paper will also shed some light on attacking the notoriously challenging problem of controlling the mdFWER under dependence for these $p$-value ordered stepwise procedures.

\section{Proofs}

\noindent {\sc Proof of Lemma \ref{lem1}.}  Let $T$ and $P$ denote the test statistic and the corresponding $p$-value for testing $H$, respectively. When testing $H$, a type 3 error occurs if $H$ is rejected and  $\theta T < 0$. Then, the type 3 error rate is given by $Pr(P \leq \alpha, \theta T < 0)$.

When $\theta >0 $, we have
\begin{eqnarray}
&& Pr(P \leq \alpha, \theta T < 0) = Pr(2F_0(T) \leq \alpha, T < 0) \nonumber \\
&=& Pr \left ( T \leq F_{0}^{-1} \left (  \frac{\alpha}{2} \right ) \right )
= F_{\theta} \left (  F_{0}^{-1} \left (  \frac{\alpha}{2} \right ) \right )\nonumber \\
&\leq& F_0 \left (  F_{0}^{-1} \left (  \frac{\alpha}{2} \right )\right )
= \frac{\alpha}{2}. \nonumber
\end{eqnarray}
The inequality follows from the assumption that $F_{\theta}$ is stochastically increasing in $\theta$. Similarly, when $\theta < 0$, we can also prove that $Pr(P \leq \alpha, \theta T < 0) \le \frac{\alpha}{2}.~\qed$

\medskip

\noindent  {\sc Proof of Theorem \ref{theorem:nodep}}(i).  Induction will be used to show
that   Procedure  1 strongly controls the mdFWER at level $\alpha$.
First consider  the case of $n=2$.  We show control of the mdFWER of Procedure \ref{procedure:nodep} in all possible combinations of true and false null hypotheses while testing two hypotheses $H_1$ and $H_2$.

\noindent \textbf{Case I: $H_{1}$ is true.} Type 1 or type 3 error occurs only when $H_1$ is rejected.
\begin{eqnarray}
\text{mdFWER} = Pr( P_{1} \le \alpha) \le  \alpha. \nonumber
\end{eqnarray}
\noindent \textbf{Case II: Both $H_{1}$ and $H_{2}$ are false.} We have no type 1 errors but only type 3 errors.
\begin{eqnarray}
\text{mdFWER} &=& Pr(\{ P_{1} \le \alpha, T_{1} \theta_{1} < 0 \} \cup \{ P_{1} \le \alpha, P_{2} \le \alpha, T_{2} \theta_{2} < 0 \}) \nonumber \\
&\leq& Pr( P_{1} \le \alpha, T_{1} \theta_{1} < 0 ) + Pr( P_{2} \le \alpha, T_{2} \theta_{2} < 0 ) \nonumber \\
& \le & \frac{\alpha}{2} + \frac{\alpha}{2} = \alpha. \nonumber
\end{eqnarray}
The first inequality follows from Bonferroni inequality and the second follows from Lemma \ref{lem1}.

\noindent \textbf{Case III: $H_{1}$ is false and $H_{2}$ is true.} The mdFWER is bounded above
by
\begin{eqnarray}
&&   Pr\text{( make type 3 error when testing }H_1) +  Pr\text{( make type 1 error when testing } H_2) \nonumber \\
& \le & Pr( P_{1} \le \alpha, T_{1} \theta_{1} < 0 ) + Pr( P_{2} \le \alpha/2 )  \nonumber \\
& \le & \frac{\alpha}{2} + \frac{\alpha}{2} = \alpha. \nonumber
\end{eqnarray}
The first inequality follows from Bonferroni inequality and the second follows from Lemma \ref{lem1} and $P_2 \sim U(0, 1)$ since $H_{2}$ is true.

Now assume  the inductive hypothesis that the mdFWER is bounded above by $\alpha$ when testing at most $n-1$ hypotheses by using Procedure \ref{procedure:nodep} at level $\alpha$. In the following, we prove the mdFWER is also bounded above by $\alpha$ when testing $n$ hypotheses $H_1, \ldots, H_n$. Without loss of generality, assume $H_1$ is a false null (if $H_1$ is a true null, the desired result directly follows by using the same argument as in Case I of $n=2$). Then, the mdFWER is bounded above by
\begin{eqnarray}
&& Pr\text{( make type 3 error when testing }H_1) \nonumber \\
&& \quad + ~Pr\text{( make at least one type 1 or type 3 errors when testing } H_2, \ldots, H_n) \nonumber \\
& \le & \frac{\alpha}{2} + \frac{\alpha}{2} = \alpha. \nonumber
\end{eqnarray}
The inequality follows from the induction assumption, noticing that $H_2, \ldots, H_n$ are tested by using Procedure \ref{procedure:nodep} at level $\alpha/2$. Thus, the desired result follows.

\noindent (ii). We now prove that the critical constants are unimprovable.
For instance, when $H_1$ is true, it is easy to see that the first critical constant, $\alpha$, is unimprovable. For each given $k= 2, \ldots, n$, when $\theta_i > 0, i=1, \ldots, k-1$ and $\theta_k = 0$, that is, $H_i, i=1, \ldots, k-1$ are false and $H_k$ is true, we present a simple joint distribution of the test statistics $T_1, \ldots, T_k$ to show that the $k$th critical constant of this procedure is also unimprovable.

Define $Z_k \sim N(0, 1)$ and $Z_i = \Phi^{-1}(|2\Phi(Z_{i+1})-1|), i = 1, \ldots, k-1$, where $\Phi(\cdot)$ is the cdf of N(0, 1). Let $q_i$ denote $Z_i$'s upper $\alpha/2^i$ quantile. It is easy to check that for each $i=1, \ldots, k$, $Z_i \sim N(0, 1)$. Thus, $-q_i$ is $Z_i$'s lower $\alpha/2^i$ quantile. In addition, by the construction of $Z_i$'s, it is easy to see that the event $Z_i \ge q_i$ is equivalent to the event $Z_{i+1} \notin (-q_{i+1}, q_{i+1})$.

Let $T_i = Z_i + \theta_i, i = 1, \ldots, k$, thus $T_i \sim N(\theta_i, 1)$. Then, as $\theta_i \rightarrow 0+$ for $i=1, \ldots, k-1$, we have
\begin{eqnarray}
\text{mdFWER} &=&\sum_{j=1}^{k-1} Pr(T_1 \ge q_1, \ldots, T_{j-1} \ge q_{j-1}, T_j \le -q_j) \nonumber \\
&& \qquad + ~Pr(T_1 \ge q_1, \ldots, T_{k-1} \ge q_{k-1}, T_k \notin (-q_k, q_k)) \nonumber \\
&=&\sum_{j=1}^{k-1} Pr(Z_1 \ge q_1, \ldots, Z_{j-1} \ge q_{j-1}, Z_j \le -q_j) \nonumber \\
&& \qquad + ~Pr(Z_1 \ge q_1, \ldots, Z_{k-1} \ge q_{k-1}, Z_k \notin (-q_k, q_k)) \nonumber
\end{eqnarray}
\begin{eqnarray}
&=&\sum_{j=1}^{k-1} Pr(Z_j \le -q_j) + Pr(Z_k \notin (-q_k, q_k)) \nonumber \\
&=&\sum_{j=1}^{k-1} \frac{\alpha}{2^j} + \frac{\alpha}{2^{(k-1)}} = \alpha. \nonumber
\end{eqnarray}
Thus, the $k$th critical constant of Procedure \ref{procedure:nodep} is unimprovable and hence each critical constant of Procedure \ref{procedure:nodep} is unimprovable under arbitrary dependence. \qed

\medskip

\noindent {\sc Proof of Lemma \ref{lem2}.} Note that when $\theta_1 >  0$ and $\theta_2 = 0$, we have
\begin{eqnarray}
&& \text{mdFWER} \nonumber \\
&=& Pr \left ( P_1 \le \alpha, \theta_1T_1 < 0 \right ) + Pr \left ( P_1 \le \alpha, \theta_1T_1 \ge 0, P_2 \le \alpha \right ) \nonumber \\
&=& Pr \left ( P_1 \le \alpha, T_1 < 0 \right ) + Pr \left ( P_1 \le \alpha, T_1 \ge 0, P_2 \le \alpha, T_2 > 0 \right ) \nonumber \\
&& + Pr \left ( P_1 \le \alpha, T_1 \ge 0, P_2 \le \alpha, T_2 \le 0 \right ) \nonumber \\
&=& Pr \left ( 2 F_0(T_1) \le \alpha \right ) + Pr \left ( 2(1-F_0(T_1)) \le \alpha, 2(1-F_0(T_2)) \le \alpha \right ) \nonumber \\
&& + Pr \left ( 2(1-F_0(T_1)) \le \alpha, 2F_0(T_2) \le \alpha \right ) \nonumber \\
&=&  Pr \left ( T_1 \leq c_1 \right ) + Pr \left ( T_1 \geq c_2, T_2 \geq c_2 \right ) + Pr \left ( T_1 \geq c_2, T_2 \leq c_1 \right ) \nonumber \\
&=& F_{\theta_1}(c_1) + 1 - F_{\theta_1}(c_2) - F_{0}(c_2) + F_{(\theta_1,0)}(c_2, c_2) + F_{0}(c_1) - F_{(\theta_1,0)}(c_2, c_1) \nonumber \\
&=& \alpha + F_{\theta_1}(c_1) - F_{\theta_1}(c_2) + F_{(\theta_1,0)}(c_2, c_2) - F_{(\theta_1,0)}(c_2, c_1).\label{equation:old10}
\end{eqnarray}
Specifically,  under Assumption \ref{assumption:indep} (independence), (\ref{equation:old10}) can be  simplified as,
\begin{eqnarray}
& & \alpha + F_{\theta_1}(c_1) - F_{\theta_1}(c_2) + F_{\theta_1}(c_2)F_{0}(c_2) - F_{\theta_1}(c_2)F_{0}(c_1) \nonumber \\
&=& \alpha + F_{\theta_1}(c_1) - \alpha F_{\theta_1}(c_2). \nonumber
\end{eqnarray}
Similarly, when $\theta_1 <  0$ and $\theta_2 = 0$, we can prove that
$$\text{mdFWER} = 1 +  F_{\theta_1}(c_1) - F_{\theta_1}(c_2) + F_{(\theta_1, 0)}(c_1, c_1) - F_{(\theta_1, 0)}(c_1, c_2). ~\qed$$

\medskip

\noindent {\sc Proof of Lemma \ref{lem3}.} By using the same arguments as in Theorem \ref{theorem:nodep}, we can easily prove control of the mdFWER of Procedure \ref{procedure:indep}
in the case of $n=2$ when $H_1$ is true or both $H_1$ and $H_2$ are false. In the following, we prove the desired result also holds when $H_{1}$ is false and $H_{2}$ is true.

Note that $H_{1}$ is false and $H_{2}$ is true imply $\theta_1 \neq 0$ and $\theta_2 = 0$. To show that the mdFWER is controlled for $\theta_1 > 0$ and $\theta_2 = 0$,  we only need to show by Lemma \ref{lem2} that $\alpha + F_{\theta_1}(c_1) - \alpha F_{\theta_1}(c_2) \le \alpha$. This is equivalent to show
\begin{eqnarray}\label{eq:old11}
F_{\theta_1}(c_2) \left ( F_0(c_2) - F_0(c_1) \right ) & \le & F_{\theta_1}(c_2) - F_{\theta_1}(c_1).
\end{eqnarray}

For proving (\ref{eq:old11}), it is enough to prove the following, as $0 \le F_{0}(c_2) \le 1$,
\begin{eqnarray}\label{eq:old12}
F_{\theta_1}(c_2) \left ( F_0(c_2) - F_0(c_1) \right )  \le F_0(c_2) \left ( F_{\theta_1}(c_2) - F_{\theta_1}(c_1) \right ).
\end{eqnarray}

Dividing both sides of (\ref{eq:old12}) by $F_{\theta_1}(c_2)F_{0}(c_2)$, we see that we only need to prove,
\begin{eqnarray}
1 - \frac{F_{0}(c_1)}{F_{0}(c_2)} & \le & 1 - \frac{F_{\theta_1}(c_1)}{F_{\theta_1}(c_2)}, \nonumber
\end{eqnarray}
which follows directly from (\ref{eq:6}) and Assumption \ref{assumption:mlr} (MLR).

Similarly, to show that the mdFWER is controlled for $\theta_1 < 0$ and $\theta_2 = 0$, we only need to show by Lemma \ref{lem2} that
$1 + \alpha F_{\theta_1}(c_1) - F_{\theta_1}(c_2) \le \alpha.$ This is equivalent to showing
\begin{eqnarray}
(1 - \alpha) \left ( 1-F_{\theta_1}(c_1) \right ) \leq F_{\theta_1}(c_2) - F_{\theta_1}(c_1). \nonumber
\end{eqnarray}
Writing $1 - \alpha$ as $\left ( 1-F_0(c_1) \right ) - \left ( 1-F_0(c_2) \right )$ and writing $ F_{\theta_1}(c_2) - F_{\theta_1}(c_1)$ as \\ $\left( 1 - F_{\theta_1}(c_1) \right) - \left ( 1-F_{\theta_1}(c_2) \right )$, we get that it is equivalent to prove
\begin{eqnarray}\label{eq:old13}
\left[ \left ( 1-F_0(c_1) \right ) - \left ( 1-F_0(c_2) \right ) \right] \left ( 1-F_{\theta_1}(c_1) \right ) \le \left( 1 - F_{\theta_1}(c_1) \right) - \left( 1-F_{\theta_1}(c_2) \right).
\end{eqnarray}
Since $0 \le 1-F_{0}(c_1) \le 1$, to prove inequality (\ref{eq:old13}), it is enough to prove the following,
\begin{eqnarray}
& & \left ( 1-F_{\theta_1}(c_1) \right ) \left[ \left ( 1-F_0(c_1) \right ) - \left ( 1-F_0(c_2) \right ) \right] \nonumber \\
& \le & \left ( 1-F_{0}(c_1) \right ) \left[ 1 - F_{\theta_1}(c_1) \right] - \left[ 1-F_{\theta_1}(c_2) \right].
\label{eq:old14}
\end{eqnarray}
Dividing both sides of (\ref{eq:old14}) by $ \left ( 1 - F_{\theta_1}(c_1) \right ) \left ( 1-F_0(c_1) \right ) $, we see that proving (\ref{eq:old13}) is equivalent to showing
\begin{eqnarray}
\frac{1-F_{\theta_1}(c_2)}{1-F_{\theta_1}(c_1)} \leq \frac{1-F_0(c_2)}{1-F_{0}(c_1)},
\end{eqnarray}
which follows directly from (\ref{eq:7}) and Assumption \ref{assumption:mlr} (MLR).
By combining the arguments of the above two cases, the desired result follows. \qed

\medskip

\noindent {\sc Proof of Theorem \ref{theorem:indep}.} The proof is by induction on number of hypotheses $n$. We already proved strong control of the mdFWER for $n=2$ in Lemma \ref{lem3}. Let us assume the result holds for testing any $n=k$ hypotheses, that is, $\text{mdFWER} \le \alpha$ while testing any $k$ pre-ordered hypotheses. We now argue that is will hold for  $n = k+1$ hypotheses. Without loss of generality, assume $H_1$ is a false null, as in the proof of Theorem \ref{theorem:nodep}.

Let $V_{k+1}^{(-1)}$ denote the total number of type 1 or type 3 errors committed while testing  $H_2, \ldots, H_{k+1}$ and excluding $H_1$. Then, by the inductive hypothesis, the mdFWER while testing the $k$ hypotheses $H_2, \ldots, H_{k+1}$ is $Pr (V_{k+1}^{(-1)} > 0) \le \alpha$. Then, the mdFWER of testing $k+1$ hypotheses $H_1, \ldots, H_{k+1}$ is defined by
\begin{eqnarray}\label{eq:old16}
& & Pr \left ( \{P_1 \le \alpha, T_1 \theta_1 < 0 \} \cup \{ P_1 \le \alpha, T_1 \theta_1 \ge 0, V_{k+1}^{(-1)} > 0  \} \right ) \nonumber \\
&=& Pr \left ( P_1 \le \alpha, T_1 \theta_1 < 0  \right ) + Pr \left ( P_1 \le \alpha, T_1 \theta_1 \ge 0 \right ) \cdot Pr \left ( V_{k+1}^{(-1)} > 0 \right )  \nonumber \\
&\le& Pr \left (P_1 \le \alpha, T_1 \theta_1 < 0 \right ) + \alpha \ Pr \left ( P_1 \le \alpha, T_1 \theta_1 \ge 0 \right ).
\end{eqnarray}
The equality follows by Assumption \ref{assumption:indep} (independence) and the inequality follows by the inductive hypothesis.  Note that (\ref{eq:old16}) is the same as (\ref{equation:old9}) under independence, which is equal to the mdFWER of Procedure \ref{procedure:indep} in the case of two hypotheses.
So again by applying Lemma \ref{lem3}, we get that $\text{mdFWER} \le \alpha$ for $n=k+1$. Hence, the proof follows by induction. ~\qed

\bigskip

\noindent {\sc Proof of  Theorem \ref{theorem:positive} .}
Without loss of generality, we assume $\theta_i > 0$ if $\theta_i \neq 0$ for $i = 1, \ldots, n$.
Also, if there exists an $i$  with  $\theta_i = 0$, by induction, we can simply assume $i_0 = n$. Thus, to prove the mdFWER control of Procedure \ref{procedure:indep}, we only need to consider two cases:

\noindent (i) $\theta_i > 0$ for $i = 1, \ldots, n$;

\noindent (ii) $\theta_i > 0$ for $i = 1, \ldots, n-1$ and
$\theta_n = 0$.

\medskip

\noindent {\bf Case (i).}
Consider the general case of $\theta_i > 0, i = 1, \ldots, n$. By Assumption \ref{assumption:prd}, the test statistics $T_1, \ldots, T_n$ are positively regression dependent. For $j=1, \ldots, n-1$, let $E_{n-j}$ denote the event of making at least one type 3 error when testing $H_{j+1}, \ldots, H_n$ using Procedure \ref{procedure:indep} at level $\alpha$. By using induction, we prove the following two lemmas hold.

\begin{lemma}\label{lem4}  Assume the conditions of Theorem \ref{theorem:positive}.
For $j=1, \ldots, n-1$, the following inequality holds.
\begin{equation}\label{eq:8}
Pr(E_{n-j}| T_1 > c_2, \ldots, T_j > c_2) \le \alpha.
\end{equation}
\end{lemma}

\noindent {\sc Proof of Lemma \ref{lem4}.} We prove the result by using reverse induction. When $j=n-1$, we have
\begin{eqnarray}
&& Pr(E_{n-j}| T_1 > c_2, \ldots, T_j > c_2) \nonumber \\
&=& Pr(T_n < c_1| T_1 > c_2, \ldots, T_{n-1} > c_2) \nonumber \\
&=& \frac{Pr(T_n < c_1)Pr(T_1 > c_2, \ldots, T_{n-1} > c_2 | T_n < c_1)}{Pr(T_1 > c_2, \ldots, T_{n-1} > c_2)} \nonumber \\
& \le & Pr(T_n < c_1) \le \alpha. \nonumber
\end{eqnarray}
The inequality follows from Assumption \ref{assumption:prd}.

Assume the inequality (\ref{eq:8}) holds for  $j=m$. In the following, we prove that it also holds for $j=m-1$. Note that
\begin{eqnarray}
&& Pr(E_{n-m+1}| T_1 > c_2, \ldots, T_{m-1} > c_2) \nonumber \\
&=& Pr \left(\{T_m < c_1\}\bigcup \left(\{T_m > c_2\} \bigcap E_{n-m}\right) \big | T_1 > c_2, \ldots, T_{m-1} > c_2 \right) \nonumber \\
&=& Pr\left(T_m < c_1 \big | T_1 > c_2, \ldots, T_{m-1} > c_2 \right) \nonumber \\
&& + ~Pr \left(\{T_m > c_2\} \bigcap E_{n-m} \big | T_1 > c_2, \ldots, T_{m-1} > c_2 \right) \nonumber \\
&=& Pr\left(T_m < c_1 \big | T_1 > c_2, \ldots, T_{m-1} > c_2 \right) \nonumber \\
&& + ~Pr \left(T_m > c_2 \big | T_1 > c_2, \ldots, T_{m-1} > c_2 \right) Pr \left(E_{n-m} \big | T_1 > c_2, \ldots, T_{m} > c_2 \right) \nonumber \\
&\le& Pr\left(T_m < c_1 \big | T_1 > c_2, \ldots, T_{m-1} > c_2 \right) + \alpha Pr \left(T_m > c_2 \big | T_1 > c_2, \ldots, T_{m-1} > c_2 \right)  \nonumber \\
&\le& \alpha. \nonumber
\end{eqnarray}
Therefore, the desired result follows. Here, the first inequality follows from the assumption of induction and the second follows from Lemma \ref{lem5} below. ~\qed

\begin{lemma}\label{lem5}Assume the conditions of Theorem \ref{theorem:positive}.
For $j=1, \ldots, n-1$, the following inequality holds:
\begin{equation}\label{eq:9}
Pr\left(T_j < c_1 \big | T_1 > c_2, \ldots, T_{j-1} > c_2 \right) + \alpha Pr \left(T_j > c_2 \big | T_1 > c_2, \ldots, T_{j-1} > c_2 \right)  \le \alpha.
\end{equation}
Specifically, for $j=1$, we have
$$Pr \left(T_1 < c_1 \right) + \alpha Pr \left(T_1 > c_2 \right) \le \alpha.$$
\end{lemma}

\noindent {\sc Proof of Lemma \ref{lem5}.} To prove the inequality (\ref{eq:9}), it is enough to show that
\begin{equation}
Pr\left(T_j < c_1 \big | T_1 > c_2, \ldots, T_{j-1} > c_2 \right) \le \alpha Pr \left(T_j < c_2 \big | T_1 > c_2, \ldots, T_{j-1} > c_2 \right), \nonumber
\end{equation}
which is equivalent to
\begin{eqnarray}\label{equation:oldstar}
&& (1 - \alpha)Pr \left(T_j < c_2 \big | T_1 > c_2, \ldots, T_{j-1} > c_2 \right) \nonumber\\
&\le& Pr \left(T_j < c_2 \big | T_1 > c_2, \ldots, T_{j-1} > c_2 \right) - Pr \left(T_j < c_1 \big | T_1 > c_2, \ldots, T_{j-1} > c_2 \right). \nonumber
\end{eqnarray}
Note that
$$1 - \alpha = Pr_{\theta_j = 0}(T_j < c_2) - Pr_{\theta_j = 0}(T_j < c_1).$$
Thus, the above inequality is equivalent to
$$
Pr_{\theta_j = 0}(T_j < c_2) - Pr_{\theta_j = 0}(T_j < c_1)  \le 1 - \frac{Pr \left(T_j < c_1 \big | T_1 > c_2, \ldots, T_{j-1} > c_2 \right)}{Pr \left(T_j < c_2 \big | T_1 > c_2, \ldots, T_{j-1} > c_2 \right)},$$
which  in turn is implied by
\begin{eqnarray}\label{eq:10}
1 - \frac{Pr_{\theta_j = 0}(T_j < c_1)}{Pr_{\theta_j = 0}(T_j < c_2)} \le 1 - \frac{Pr \left(T_j < c_1 \big | T_1 > c_2, \ldots, T_{j-1} > c_2 \right)}{Pr \left(T_j < c_2 \big | T_1 > c_2, \ldots, T_{j-1} > c_2 \right)}.
\end{eqnarray}
Note that by Assumption \ref{assumption:mlr}, we have
$$\frac{Pr(T_j < c_1)}{Pr(T_j < c_2)} \le \frac{Pr_{\theta_j = 0}(T_j < c_1)}{Pr_{\theta_j = 0}(T_j < c_2)}.$$
Thus, to prove the inequality (\ref{eq:10}), we only need to show that
$$\frac{Pr \left(T_j < c_1 \big | T_1 > c_2, \ldots, T_{j-1} > c_2 \right)}{Pr \left(T_j < c_2 \big | T_1 > c_2, \ldots, T_{j-1} > c_2 \right)} \le \frac{Pr(T_j < c_1)}{Pr(T_j < c_2)},$$
which is equivalent to
$$Pr \left(T_1 > c_2, \ldots, T_{j-1} > c_2 \big | T_j < c_1 \right) \le Pr \left(T_1 > c_2, \ldots, T_{j-1} > c_2 \big | T_j < c_2 \right),$$
which follows from Assumption \ref{assumption:prd}. Therefore, the desired result follows. ~\qed

\vskip 5pt

Based on Lemmas \ref{lem4} and \ref{lem5}, we have
\begin{eqnarray}
&& \text{mdFWER}  = Pr(T_1 < c_1) + \sum_{j=2}^{n} Pr(T_1 > c_2, \ldots, T_{j-1} > c_2, T_j < c_1) \nonumber \\
&=& Pr(T_1 < c_1) + Pr(T_1 > c_2) \sum_{j=2}^{n} Pr(T_2 > c_2, \ldots, T_{j-1} > c_2, T_j < c_1| T_1 > c_2) \nonumber \\
&=& Pr(T_1 < c_1) + Pr(T_1 > c_2) Pr(E_{n-1}| T_1 > c_2) \nonumber \\
&\le& Pr(T_1 < c_1) + \alpha Pr(T_1 > c_2) \nonumber \\
&\le& \alpha. \nonumber
\end{eqnarray}
Therefore, the mdFWER is controlled at level $\alpha$ for Case (i).
Here, the first inequality follows from Lemma \ref{lem4} and the second follows from Lemma \ref{lem5}.

\medskip

\noindent{\bf Case (ii).}  Consider the general case of $\theta_i >0, i=1, \ldots, n-1$ and $\theta_n = 0$. Under Assumption \ref{assumption:prd}, $T_i, i=1, \ldots, n-1$ are positively regression dependent and under Assumption \ref{assumption:indep2}, $T_n$ is independent of $T_i$'s . Note that
\begin{eqnarray}
&& \text{mdFWER}  \nonumber \\
&=& \sum_{j=1}^{n-1} Pr(T_1 > c_2, \ldots, T_{j-1} > c_2, T_j < c_1) \nonumber \\
&& \quad + ~Pr(T_1 > c_2, \ldots, T_{n-1} > c_2, T_n < c_1) + ~Pr(T_1 > c_2, \ldots, T_n > c_2) \nonumber \\
&=& \sum_{j=1}^{n-1} Pr(T_1 > c_2, \ldots, T_{j-1} > c_2, T_j < c_1) + \alpha Pr(T_1 > c_2, \ldots, T_{n-1} > c_2). \nonumber
\end{eqnarray}
The second equality follows from Assumption \ref{assumption:indep2}. \\
For $m=1, \ldots, n-1$, define
$$\Delta_m = \sum_{j=1}^{m} Pr(T_1 > c_2, \ldots, T_{j-1} > c_2, T_j < c_1) + \alpha Pr(T_1 > c_2, \ldots, T_m > c_2).$$
Thus, $\text{mdFWER} = \Delta_{n-1}$. By using induction, we prove below that $\Delta_m \le \alpha$ for $m=1, \ldots, n-1$.

For $m=1$, by using Lemma \ref{lem5}, we have
$$\Delta_1 = Pr \left(T_1 < c_1 \right) + \alpha Pr \left(T_1 > c_2 \right) \le \alpha.$$
Assume $\Delta_m \le \alpha$. In  the following, we show $\Delta_{m+1} \le \alpha$. Note that
\begin{eqnarray}
\Delta_{m+1} &=& \sum_{j=1}^{m+1} Pr(T_1 > c_2, \ldots, T_{j-1} > c_2, T_j < c_1) \nonumber \\
&& \quad +~ \alpha Pr(T_1 > c_2, \ldots, T_m > c_2, T_{m+1} > c_2) \nonumber \\
&=& \sum_{j=1}^{m} Pr(T_1 > c_2, \ldots, T_{j-1} > c_2, T_j < c_1) \nonumber \\
&& \quad +~ Pr(T_1 > c_2, \ldots, T_m > c_2) \left [Pr(T_{m+1} < c_1|T_1 > c_2, \ldots, T_m > c_2) \right .  \nonumber \\
&& \qquad \left . +~  \alpha Pr(T_{m+1} > c_2|T_1 > c_2, \ldots, T_m > c_2) \right ] \nonumber \\
&\le& \sum_{j=1}^{m} Pr(T_1 > c_2, \ldots, T_{j-1} > c_2, T_j < c_1) + \alpha Pr(T_1 > c_2, \ldots, T_m > c_2) \nonumber \\
&=& \Delta_{m} \le \alpha.
\end{eqnarray}
The first inequality follows from Lemma \ref{lem5} and the second follows from the inductive hypothesis. Thus, $\Delta_{m} \le \alpha$ for $m=1, \ldots, n-1$. Therefore, $\text{mdFWER} = \Delta_{n-1} \le \alpha$, the desired result.

\vskip 5pt

Combining the arguments  of Cases (i) and (ii), the proof of Theorem \ref{theorem:positive} is complete. ~\qed

\bigskip

\noindent {\sc Proof of Proposition \ref{proposition:1}}. From the proof of Theorem \ref{theorem:nodep} and by Lemma \ref{lem1}, it is easy to see that we only need to prove the mdFWER control of Procedure \ref{procedure:indep} when $H_1$ is false and $H_2$ is true, i.e., $\theta_1 \neq 0$ and $\theta_2 = 0$.

\noindent \textbf{Case I: $\theta_1 > 0$ and $\theta_2 = 0$.} By Lemma \ref{lem2}, the mdFWER of Procedure \ref{procedure:indep} is controlled at level $\alpha$ if we have the following:
\begin{eqnarray}
F_{\theta_1}(c_1) - F_{\theta_1}(c_2) + F_{(\theta_1,0)}(c_2,c_2) - F_{(\theta_1,0)}(c_2,c_1) \le 0. \nonumber
\end{eqnarray}
After rewriting $F_{(\theta_1,0)}(x,y)$ as $Pr(T_1\le x, T_2\le y)$ and then dividing  through by $Pr(T_1 \le c_2)$, we get,
\begin{eqnarray}
Pr \left ( T_2 \leq c_2 | T_1 \leq c_2 \right ) - Pr \left ( T_2 \leq c_1 | T_1 \leq c_2 \right ) \leq 1- \frac{Pr(T_1 \leq c_1)}{Pr(T_1 \leq c_2)}. \nonumber
\end{eqnarray}
Dividing  by $Pr \left ( T_2 \leq c_2 | T_1 \leq c_2 \right )$, we get,
\begin{eqnarray}\label{eq:11}
1- \frac{Pr \left ( T_2 \leq c_1 | T_1 \leq c_2 \right )}{Pr \left ( T_2 \leq c_2 | T_1 \leq c_2 \right )} \leq \frac{1}{Pr \left ( T_2 \leq c_2 | T_1 \leq c_2 \right )} \left ( 1- \frac{Pr(T_1 \leq c_1)}{Pr(T_1 \leq c_2)} \right ).
\end{eqnarray}
For proving (\ref{eq:11}), it is enough to prove the following inequality, as $\frac{1}{Pr \left ( T_2 \le c_2 | T_1 \le c_2 \right )} \ge 1 $.
\begin{eqnarray}\label{eq:12}
1- \frac{Pr \left( T_2 \le c_1 | T_1 \le c_2 \right )}{Pr \left( T_2 \le c_2 | T_1 \le c_2 \right )} \le 1- \frac{Pr(T_1 \le c_1)}{Pr(T_1 \le c_2)}.
\end{eqnarray}
By Assumption \ref{assumption:mlr} and (\ref{eq:6}), it follows that $\frac{F_{0}(c_2)}{F_{0}(c_1)} \le \frac{F_{\theta_1}(c_2)}{F_{\theta_1}(c_1)}$, which is equivalent to, $1-\frac{Pr(T_2 \le c_1)}{Pr(T_2\leq c_2)} \leq 1-\frac{Pr(T_1 \leq c_1)}{Pr(T_1 \leq c_2)}$. Thus for proving (\ref{eq:11}), it is enough to prove the following:
\begin{eqnarray}\label{equation:old43}
1- \frac{Pr \left( T_2 \le c_1 | T_1 \le c_2 \right )}{Pr \left( T_2 \le c_2 | T_1 \le c_2 \right )} \le 1- \frac{Pr(T_2 \le c_1)}{Pr(T_2 \le c_2)}~.
\end{eqnarray}
But, (\ref{equation:old43})  is  equivalent to showing
\begin{eqnarray}
Pr \left( T_1 \le c_2 | T_2 \le c_1 \right ) \ge Pr \left( T_1 \le c_2 | T_2 \le c_2 \right ), \nonumber
\end{eqnarray}
which follows directly from Assumption \ref{assumption:prds2}.

\noindent \textbf{Case II: $\theta_1 < 0$ and $\theta_2 = 0$.} Similarly, by Lemma \ref{lem2}, the mdFWER of Procedure \ref{procedure:indep} is controlled at level $\alpha$ if we have the following:
\begin{eqnarray}\label{eq:12}
1 +  F_{\theta_1}(c_1) - F_{\theta_1}(c_2) + F_{(\theta_1, 0)}(c_1, c_1) - F_{(\theta_1, 0)}(c_1, c_2) \le \alpha,
\end{eqnarray}
which after some rearrangement and rewriting $1-\alpha$ as $F_{0}(c_2) - F_{0}(c_1)$ gives,
\begin{eqnarray}
\left( F_{0}(c_2) - F_{(\theta_1, 0)}(c_1, c_2) \right ) - \left( F_{0}(c_1) - F_{(\theta_1, 0)}(c_1, c_1) \right ) \le \left( 1 -  F_{\theta_1}(c_1) \right ) - \left( 1 - F_{\theta_1}(c_2) \right ).
\end{eqnarray}
Thus, proving (\ref{eq:12}) is equivalent to proving that
\begin{eqnarray}
Pr \left( T_1 \ge c_1, T_2 \le c_2 \right ) - Pr \left( T_1 \ge c_1, T_2 \le c_1 \right ) \le Pr \left( T_1 \ge c_1 \right ) - Pr \left( T_1 \ge c_2 \right ). \nonumber
\end{eqnarray}
Dividing through by $Pr(T_1 \ge c_1)$, we get
\begin{eqnarray}
Pr \left( T_2 \ge c_1 | T_1 \ge c_1 \right ) - Pr \left( T_2 \ge c_2 | T_1 \ge c_1 \right ) \le 1 - \frac{Pr(T_1 \ge c_2)}{Pr(T_1 \ge c_1)}.
\end{eqnarray}
Thus to prove (\ref{eq:12}), it is enough to prove the following,
\begin{eqnarray}
1 - \frac{Pr \left( T_2 \ge c_2 | T_1 \ge c_1 \right )}{Pr \left( T_2 \ge c_1 | T_1 \ge c_1 \right )} \le 1 - \frac{Pr(T_1 \ge c_2)}{Pr(T_1 \ge c_1)}, \nonumber
\end{eqnarray}
which is equivalent to proving,
\begin{eqnarray}
\frac{Pr \left( T_2 \ge c_2 | T_1 \ge c_1 \right )}{Pr \left( T_2 \ge c_1 | T_1 \ge c_1 \right )} \ge \frac{Pr(T_1 \ge c_2)}{Pr(T_1 \ge c_1)}.
\end{eqnarray}
By Assumption \ref{assumption:mlr} and (\ref{eq:7}), it follows that for $\theta_1 < 0$, $\frac{Pr(T_1 \ge c_2)}{Pr(T_1 \ge c_1)} \le \frac{Pr(T_2 \ge c_2)}{Pr(T_2 \ge c_1)}$. Thus to prove (\ref{eq:12}), it is enough to prove the following,
\begin{eqnarray}
\frac{Pr \left( T_2 \ge c_2 | T_1 \ge c_1 \right )}{Pr \left( T_2 \ge c_1 | T_1 \ge c_1 \right )} \ge \frac{Pr(T_2 \ge c_2)}{Pr(T_2 \ge c_1)}. \label{equation:newnumber1}
\end{eqnarray}
But (\ref{equation:newnumber1})  is equivalent to showing
\begin{eqnarray}
Pr \left( T_1 \ge c_1 | T_2 \ge c_2 \right ) \ge Pr \left( T_1 \ge c_1 | T_2 \ge c_1 \right ),
\end{eqnarray}
which follows directly from Assumption \ref{assumption:prds2}.
By combining the arguments of the above two cases, the desired result follows. ~\qed

\bigskip
\noindent{\sc Proof of Proposition \ref{proposition:2}}.
By Corollary \ref{corollary:1}, without loss of generality, assume that $\theta_i > 0, i=1, 2$ and $\theta_3 =0$, that is, $H_1$ and $H_2$ are false and $H_3$ is true. Note that
\begin{eqnarray}\label{eq:12_2}
&& \text{mdFWER}  \\
&=& Pr(T_1 \le c_1) + Pr(T_1 \ge c_2, T_2 \le c_1) + Pr\left (T_1 \ge c_2, T_2 \ge c_2, T_3 \notin (c_1, c_2) \right). \nonumber
\end{eqnarray}
In the following, we prove that
\begin{eqnarray}\label{eq:13}
&& Pr(T_1 \ge c_2, T_2 \le c_1) + Pr\left (T_1 \ge c_2, T_2 \ge c_2, T_3 \notin (c_1, c_2) \right) \nonumber \\
& \le & Pr\left (T_1 \ge c_2, T_3 \notin (c_1, c_2) \right).
\end{eqnarray}
To prove (\ref{eq:13}), it is enough to show the following inequality:
\begin{eqnarray}\label{eq:14}
Pr(T_2 \le c_1|T_1) + Pr\left (T_2 \ge c_2, T_3 \notin (c_1, c_2) |T_1 \right) \le Pr\left (T_3 \notin (c_1, c_2) |T_1 \right).
\end{eqnarray}
Note that
\begin{eqnarray}\label{eq:15}
Pr\left (T_2 \ge c_2, T_3 \le c_1 |T_1 \right) =  Pr(T_3 \le c_1|T_1) - Pr\left (T_2 < c_2, T_3 \le c_1 |T_1\right)
\end{eqnarray}
and
\begin{eqnarray}\label{eq:16}
&& Pr\left (T_2 \ge c_2, T_3 \ge c_2 |T_1 \right) \nonumber \\
& = & 1 - Pr(T_2 < c_2 |T_1) - Pr(T_3 < c_2 |T_1) + Pr\left (T_2 < c_2, T_3 < c_2 |T_1 \right).
\end{eqnarray}
In addition, we have
\begin{eqnarray}\label{eq:17}
Pr\left (T_3 \notin (c_1, c_2) |T_1 \right) = 1 + Pr(T_3 \le c_1|T_1) - Pr(T_3 < c_2|T_1).
\end{eqnarray}
Thus, in order to show (\ref{eq:14}), by combining (\ref{eq:15})-(\ref{eq:17}), we only need to prove the following inequality:
\begin{eqnarray}\label{eq:18}
&& Pr\left (T_2 < c_2, T_3 < c_2 |T_1 \right) - Pr\left (T_2 < c_2, T_3 \le c_1 |T_1 \right) \nonumber \\
& \le & Pr(T_2 < c_2 |T_1) - Pr(T_2 \le c_1 |T_1).
\end{eqnarray}
Note that (\ref{eq:18}) can be rewritten as
\begin{eqnarray}\label{eq:19}
&& Pr\left (T_2 < c_2, T_3 < c_2 |T_1 \right) \left [1  - \frac{Pr\left (T_2 < c_2, T_3 \le c_1 |T_1 \right)}{Pr\left (T_2 < c_2, T_3 < c_2 |T_1 \right)} \right ] \nonumber \\
& \le & Pr(T_2 < c_2 |T_1) \left[ 1 - \frac{Pr(T_2 \le c_1 |T_1)}{Pr(T_2 < c_2 |T_1)} \right].
\end{eqnarray}
Thus, to prove (\ref{eq:18}), it is enough to show
\begin{eqnarray}\label{eq:20}
1  - \frac{Pr\left (T_2 < c_2, T_3 \le c_1 |T_1 \right)}{Pr\left (T_2 < c_2, T_3 < c_2 |T_1 \right)} \le 1 - \frac{Pr(T_2 \le c_1 |T_1)}{Pr(T_2 < c_2 |T_1)}.
\end{eqnarray}
That is,
\begin{eqnarray}\label{eq:21}
\frac{Pr(T_2 \le c_1 |T_1)}{Pr(T_2 < c_2 |T_1)} \le \frac{Pr\left (T_2 < c_2, T_3 \le c_1 |T_1 \right)}{Pr\left (T_2 < c_2, T_3 < c_2 |T_1 \right)}.
\end{eqnarray}

By Assumption \ref{assumption:bmlr} (BMLR), we have
\begin{eqnarray}\label{eq:24}
\frac{Pr(T_2 \le x_2|T_1)}{Pr(T_3 \le x_2|T_1)} \ge \frac{Pr(T_2 \le x_1|T_1)}{Pr(T_3 \le x_1|T_1)}.
\end{eqnarray}
By (\ref{eq:24}), to prove (\ref{eq:21}), it is enough to show
\begin{eqnarray}\label{eq:25}
\frac{Pr(T_3 \le c_1 |T_1)}{Pr(T_3 < c_2 |T_1)} \le \frac{Pr\left (T_2 < c_2, T_3 \le c_1 |T_1 \right)}{Pr\left (T_2 < c_2, T_3 < c_2 |T_1 \right)}.
\end{eqnarray}
That is,
\begin{eqnarray}\label{eq:26}
Pr\left (T_2 < c_2 |T_3 < c_2, T_1 \right) \le Pr\left (T_2 < c_2 |T_3 < c_1, T_1 \right).
\end{eqnarray}
The inequality (\ref{eq:26}) holds under Assumption \ref{assumption:prds2}. Therefore, the inequality (\ref{eq:13}) holds.

Based on (\ref{eq:12_2})-(\ref{eq:13}) and Proposition 1, we have
\begin{eqnarray}
\text{mdFWER} = Pr(T_1 \le c_1) + Pr\left (T_1 \ge c_2, T_3 \notin (c_1, c_2) \right) \le \alpha. \nonumber
\end{eqnarray}
Thus, the desired result follows. ~\qed
 \bigskip

\noindent{\sc Proof of Theorem \ref{theorem:n}}.
By Corollary \ref{corollary:1}, without loss of generality, assume that $\theta_i > 0, i=1, \ldots, n-1$ and $\theta_n =0$, that is, $H_i, i=1, \ldots, n-1$ are false and $H_n$ is true. Note that
\begin{eqnarray}\label{eq:27}
&& \text{mdFWER}   \\
&=& \sum_{j=1}^{n-1} Pr(T_1 \ge c_2, \ldots, T_{j-1} \ge c_2, T_j \le c_1) + Pr(T_1 \ge c_2, \ldots, T_{n-1} \ge c_2, T_n \notin (c_1, c_2)). \nonumber
\end{eqnarray}
In the following, we prove that
\begin{eqnarray}\label{eq:28}
&& Pr(T_1 \ge c_2, \ldots, T_{n-2} \ge c_2, T_{n-1} \le c_1) + Pr\left (T_1 \ge c_2, \ldots, T_{n-1} \ge c_2, T_n \notin (c_1, c_2) \right) \nonumber \\
& \le & Pr\left (T_1 \ge c_2, \ldots, T_{n-2} \ge c_2, T_n \notin (c_1, c_2) \right).
\end{eqnarray}
To prove (\ref{eq:28}), it is enough to show the following inequality:
\begin{eqnarray}\label{eq:29}
&& Pr(T_{n-1} \le c_1|T_1, \ldots, T_{n-2}) + Pr\left (T_{n-1} \ge c_2, T_n \notin (c_1, c_2) |T_1, \ldots, T_{n-2} \right) \nonumber \\
&\le& Pr\left (T_n \notin (c_1, c_2) |T_1, \ldots, T_{n-2} \right).
\end{eqnarray}
By using the same argument as in proving (\ref{eq:14}) in the case of three hypotheses, we can prove that the inequality (\ref{eq:29}) holds under Assumptions \ref{assumption:prds2} and \ref{assumption:mmlr}. Then, by combining (\ref{eq:27}) and (\ref{eq:28}), we have
\begin{eqnarray}\label{eq:31}
&& \text{mdFWER}   \\
&\le& \sum_{j=1}^{n-2} Pr(T_1 \ge c_2, \ldots, T_{j-1} \ge c_2, T_j \le c_1) + Pr(T_1 \ge c_2, \ldots, T_{n-2} \ge c_2, T_n \notin (c_1, c_2)). \nonumber
\end{eqnarray}
Note that the right-hand side of (\ref{eq:31}) is the mdFWER of Procedure \ref{procedure:indep} when testing $H_1, \ldots, H_{n-2}$, $H_n$. By induction and Proposition 1, the mdFWER is bounded above by $\alpha$, the desired result.

\vskip 10pt

\section*{Acknowledgements}
The research of Wenge Guo was supported in part
by NSF Grant DMS-1309162 and the
research of Joseph Romano was supported in part
by NSF Grant DMS-0707085. We sincerely thank a referee for giving helpful and
insightful comments and Yalin Zhu for implementing the proposed procedures
in the R package FixSeqMTP.

\end{document}